\documentstyle[leqno,twoside,12pt]{article}
\newcommand{\proof}{{\em Proof.\ }}

\newcommand{\ra}{\rightarrow}

\newcommand{\arr}[1]{\stackrel{#1}{\rightarrow}}

\newcommand{\larray}{\left(\begin{array}{cc}}
\newcommand{\rarray}{\end{array}\right)}

\newcommand{\beq}{\begin{equation}}
\newcommand{\eeq}{\end{equation}}
\newcommand{\barr}{\begin{array}}
\newcommand{\earr}{\end{array}}
\newcommand{\beqar}{\begin{eqnarray}}
\newcommand{\eeqar}{\end{eqnarray}}
\newtheorem{theorem}{Theorem}

\newtheorem{lemma}[theorem]{Lemma}

\newcommand{\bc}{{\bf C}}

\newcommand{\bt}{{\bf T}}

\newcommand{\frechet}{Fr\'{e}chet}
\newcommand{\Tot}{\mbox{\rm Tot}}

\newcommand{\tdb}{\tilde{b}}
\newcommand{\tdB}{\tilde{B}}
\newcommand{\dlim}{\displaystyle{\lim_{\longrightarrow}}\,}

\newcommand{\hhci}[2]{HH_{#1}(#2, \bar{\otimes})}
\newcommand{\hcci}[2]{HC_{#1}(#2, \bar{\otimes})}
\newcommand{\hpci}[2]{HP_{#1}(#2, \bar{\otimes})}

\begin{document}

\begin{flushleft}
{\Large Periodic cyclic homology of certain nuclear algebras}\\
\mbox{}\\
{ Jacek BRODZKI \raisebox{1ex}{\scriptsize a},   Roger PLYMEN
 \raisebox{1ex}{\scriptsize
b}}\\
\rule{0mm}{2.5mm}\\
\scriptsize \raisebox{.6ex}{\scriptsize a} 
School of Mathematical Sciences, University of Exeter, North Park Road, Exeter, EX4 4QE, U.K
\newline E-mail: brodzki@maths.ex.ac.uk\\
\raisebox{.6ex}{\scriptsize b} Department of Mathematics, University of Manchester,
Manchester, M13 9PL, U.K.,
\newline E-mail: roger@ma.man.ac.uk
\end{flushleft}

\vspace{0.3cm}

\begin{center}
\parbox[t]{2cm}{\bf Abstract}\\

\mbox{}

\parbox[t]{11.5cm}{\footnotesize
Relying on  properties of the inductive tensor product, 
we construct cyclic type homology theories for certain nuclear algebras.  
In this context, we establish continuity theorems.
 We compute the periodic cyclic homology of the 
Schwartz algebra of $p$-adic $GL(n)$ in terms of compactly supported 
de Rham cohomology of the tempered dual of $GL(n)$.}

\end{center}

\noindent{\bf 1. Complete nuclear locally convex algebras} 

Cyclic type homology groups of an algebra $A$ are computed 
using chain complexes involving tensor powers of $A$. When $A$ is a general 
locally 
convex algebra, this will involve making a choice of a topological tensor 
product. 
A {\em locally convex algebra}  is a locally
convex vector space
$A$ over $\bc$ equipped with a separately continuous multiplication. 
 We shall refer to the {\em
projective} tensor product $\otimes_\pi$, the {\em injective} tensor product
$\otimes_\epsilon$, and the {\em inductive} tensor product $\otimes_i$. Let $E$ denote a
locally convex space. If $E$ is nuclear then 
$E\otimes_\pi F \simeq E\otimes_\epsilon F$. This is the defining property of nuclear spaces
\cite[II.34]{Grot}.  The projective tensor product solves the universal problem for continuous
bilinear maps; the inductive tensor product solves the universal problem for {\em
separately} continuous bilinear maps.
If $E$ is a \frechet\ space, then $E\otimes_\pi E \simeq
E\otimes_i E$ by \cite[III.30, Corollary 1]{Bou:TVS}.  

Let $M$ be a compact $C^\infty$-manifold and let $E = C^\infty (M)$ furnished
with its standard seminorm topology. Then $E$ is nuclear and \frechet.
Therefore, the class of topologies compatible, in the sense of Grothendieck 
\cite[I.89]{Grot}, with the tensor product structure on $E\otimes E$, is a class 
with one element.  It is with respect to this unique topological 
tensor product that the cyclic homology of
the locally convex unital algebra $C^\infty(M)$ was computed by Connes \cite[Ch.II,
Theorem 46]{Connes}. 

Let $S(G)$ be the Schwartz algebra of a reductive $p$-adic group. Then 
$
S(G) = \bigcup_K S(G//K)
$
in the inductive limit topology, where $K$ is a compact open subgroup of $G$. 
The space $S(G)$ is a complete Hausdorff nuclear
topological vector space equipped with a separately continuous multiplication. 
It is the strict inductive limit of unital nuclear \frechet\ algebras $S(G//K)$. 
When we turn to the cyclic homology of $S(G)$, we are faced with a choice of
topological tensor product. Topological tensor products for nuclear
spaces such as $S(G)$ are {\em not} unique \cite[II.85]{Grot}.
We choose the completed inductive tensor product 
$\bar{\otimes}$, 
as this has good compatibility with strict inductive limits 
\cite[I.76, Prop.~14]{Grot}. This  compatibility is used in a
crucial way throughout this Note.

\begin{theorem}\label{nuke}
Let $A$ be the strict inductive limit of the nuclear \frechet\ algebras $A_\alpha$ 
with $\alpha = 1,2, 3, \dots$. Then 

{\rm (1)} $A$ is a complete Hausdorff nuclear locally
convex algebra

{\rm (2)} For all $n\geq 1$, $A^{\bar{\otimes}n}$ is 
a complete Hausdorff nuclear locally convex space and $A^{\bar{\otimes } n}= \dlim
(A_\alpha^{\bar{\otimes}n})$.
\end{theorem}
\proof (1) We may suppose that the $A_\alpha$ form an increasing sequence of 
vector subspaces of $A$ such that $A= \bigcup A_\alpha$ as in \cite[I.12, I.13]{Grot}. 
 Now $A_\alpha$ is closed in $A_{\alpha +1}$ by 
definition \cite[I.12]{Grot} so $A_\alpha$ is closed in $A$ \cite[II.32,
Prop.~9]{Bou:TVS}. Let $y_n\to 0$ in $A$, then $(y_n)$ is a bounded set hence there 
exists  an $\alpha$ for which $y_n\in A_\alpha$ by \cite[III.5, Prop.~6]{Bou:TVS}. 
Let $a\in A$ then $a\in A_\beta$ so take $\gamma = \max (\alpha , \beta)$. 
Then $a, y_n\in A_\gamma$ and $y_n\to 0$ in $A_\gamma$. So $ay_n\to 0$ in $A_\gamma$
by separate continuity of multiplication in $A_\gamma$. Then $ay_n\to 0$ in 
$A$, so $A$ is a locally convex algebra. Also $A$ is complete and Hausdorff by \cite[II.32,
Prop.~9]{Bou:TVS} and nuclear  by \cite[II.48, Corollaire 1]{Grot}.

(2) Since $A_\alpha$ is \frechet\ we have $A_\alpha \bar{\otimes} A_\alpha
= A_\alpha \hat{\otimes} A_\alpha$ by \cite[I.74]{Grot}, where $\hat{\otimes}$ 
is the completed projective tensor product. Then
$A_\alpha\bar{\otimes} A_\alpha$ is nuclear \cite[II.47, Th\'{e}or\`{e}me 9]{Grot} and
\frechet\ \cite[I.43,  Prop.~5]{Grot}. The Collorary in  \cite[II.70]{Grot} implies 
that $(A_\alpha \bar{\otimes} A_\alpha)$ is a strict inductive system. Then $\lim
(A_\alpha \bar{\otimes} A_\alpha)$ is complete by \cite[II.32, Prop.~9]{Bou:TVS} 
and so $A\bar{\otimes} A = \lim (A_\alpha \bar{\otimes} A_\alpha)$ by 
\cite[I.76, Prop.~14]{Grot}. Then 
$A\bar{\otimes} A$ is the strict inductive limit of nuclear \frechet\ spaces
hence is a complete Hausdorff nuclear locally convex space, as in  (1). 

 An argument on similar lines shows that $A^{\bar{\otimes} n} = \lim
A_\alpha^{\bar{\otimes} n}$ for all $n\geq 1$. Then $A^{\bar{\otimes}n }$ 
is the strict inductive limit of nuclear \frechet\ spaces hence is  a complete
Hausdorff nuclear locally convex space, as in (1).

\mbox{}

{\bf 2. Cyclic homology} 

Let $A$ be a locally convex algebra; we do not assume that $A$ has a unit. We denote by
$\tilde{A} = \bc \oplus A$ the unitization  of $A$. 
We associate with $A$ the mixed complex $(\bar{\Omega}\tilde{A}, \tilde{b}, \tilde{B})$ 
of noncommutative differential forms \cite{CQ2}, see also \cite{brodzki}. In positive
degrees,
$\bar{\Omega}^n\tilde{A} = A^{\bar{\otimes} n+1} \oplus A^{\bar{\otimes} n}$. 
We put $\bar{\Omega}^0\tilde{A} = A$ 
and $\bar{\Omega}^n\tilde{A} = 0$ for negative $n$. 
The differentials $\tdb$ and $\tdB$, of degree $-1$ and $+1$, respectively, 
are given  by
$$
\tilde{b} = \larray b & 1-\lambda \\ 0  & -b'\rarray ,
\qquad 
\tilde{B} = \larray 0 & 0 \\ N_\lambda & 0 \rarray. 
$$
The continuous  differentials $b'$ and $b$ of degree $-1$ are, for  $n >0$,  given by 
$$ \begin{array}{rcl}
b'(a_1\otimes \cdots \otimes a_n) & = &  \sum_{i=1}^{n-1}(-1)^{i+1}a_1\otimes \cdots \otimes a_i a_{i+1}\otimes
\cdots \otimes a_n,\\
\rule{0mm}{4mm}
b (a_1\otimes \cdots \otimes a_n) & = & b'(a_1\otimes \cdots \otimes a_n) + (-1)^{n-1}
a_na_1\otimes \cdots \otimes a_{n-1}.
\end{array}
$$
Since the (signed) generator $\lambda $ of cyclic permutations of $A^{\bar{\otimes}n}$ is 
 continuous, then so is the operator $N_\lambda = \sum_{i=0}^{n-1} \lambda^i$. Thus the 
differentials $\tdb$ and $\tdB$ are continuous.
Moreover, we have that $\tdb^2 = \tdb \tdB + \tdB \tdb = \tdB^2 = 0$.

When defining cyclic type homology theories we shall, {\em unless the topological 
tensor
product is unique}, 
 indicate explicitly the 
topological tensor product used. 
{\em Hochschild homology} $\hhci{*}{A}$ of the algebra $A$, computed with 
respect to  $\bar{\otimes}$, is by definition 
the homology of the complex
$(\bar{\Omega}\tilde{A}, \tdb)$. {\em Cyclic homology}  is defined as
$
\hcci{*}{A} = H_*(\Tot \bar{\Omega}\tilde{A} , \tdb + \tdB),$ where 
$\Tot \bar{\Omega}\tilde{A}$ is 
the total complex of the double complex associated with the
mixed complex $(\bar{\Omega}\tilde{A}, \tdb, \tdB)$. In degree $n$ we have the finite 
sum (in the direct sum topology)
$$\Tot_n\bar{\Omega}\tilde{A} = \bigoplus_{p\geq 0}\bar{\Omega}^{n-2p}
\tilde{A}.
$$
The differential $\tdb +\tdB$ is continuous
in this topology. 
 Finally, the {\em periodic cyclic homology}
 $\hpci{*}{A}$ of $A$ is the homology of the complex
$$
\cdots\arr{\tdb + \tdB} \bar{\Omega}^{even}\tilde{A}\arr{\tdb + \tdB}
\bar{\Omega}^{odd}\tilde{A}
\arr{\tdb + \tdB} \bar{\Omega}^{even}\tilde{A} \arr{\tdb + \tdB}\cdots
$$
where the spaces of even/odd chains 
$$
\bar{\Omega}^{even}\tilde{A} = \prod_{n\geq 0}\bar{\Omega}^{2n}\tilde{A}, 
\qquad
\bar{\Omega}^{odd}\tilde{A} = \prod_{n\geq 0}\bar{\Omega}^{2n+1}\tilde{A} 
$$
are equipped with the product topology which makes the differential $\tdb + \tdB$ continuous. 
\begin{theorem}\label{HHcts}
Let $A$ and $ A_\alpha$ be as in Theorem \ref{nuke}.  Then 
$$
\begin{array}{rcl}
\hhci{*}{A} & = & \displaystyle{\dlim HH_*(A_\alpha)}
 \\
HC_*( A, \bar{\otimes}) & = & 
\displaystyle{\dlim HC_*(A_\alpha)}.
\end{array}
$$
\end{theorem}
\proof Theorem \ref{nuke} gives that there is a strict inductive system of mixed complexes
$(\bar{\Omega}\tilde{A}_\alpha, \tdb_\alpha, \tdB_\alpha)$ such that, for any $n\geq 0$,
$
\displaystyle{\dlim} \bar{\Omega}^n\tilde{A}_\alpha =
 \bar{\Omega}^n(\displaystyle{\dlim} \tilde{A}_\alpha)
= \bar{\Omega}^n \tilde{A}.
$
Using the fact that homology commutes with direct limits \cite[p.~28, Prop.~1]{Bou:Ch10} and 
this remark we have
$$
\dlim HH_*(A_\alpha) = \dlim H_*(\bar{\Omega}\tilde{A}_\alpha,
\tdb_\alpha) = H_*(\dlim\bar{\Omega}
\tilde{A}_\alpha, \dlim \tdb_\alpha) = 
H_*(\bar{\Omega}\tilde{A}, \tdb) = 
HH_*(A, \bar{\otimes})
$$
where $\tdb = \dlim\tdb_\alpha$. This differential is continuous by 
Theorem \ref{nuke}.  Continuity of 
cyclic homology is proved in the same way when we use the fact that 
direct limits commute with direct sums. 

\begin{theorem}\label{HPcts}
Let $A$ and $ A_\alpha$ be as in Theorem \ref{nuke}. 
Assume that there exists $N>0$ such that $HH_n(A_\alpha) = 0$
for all $n> N$ and all $\alpha$. Then 
$$
HP_*(A, \bar{\otimes}) = \dlim  HP_*(A_\alpha).
$$
\end{theorem}
\proof Let $D$ be a locally convex algebra such that $HH_n(D) = 0$ 
for all $n > N$. Then $HP_{even}(D) = HC_{2n}(D)$, $HP_{odd}(D)
= HC_{2n+1}(D)$ for any $n$ such that $2n$ and $2n+1$ are greater than $N$. 
Indeed, let us define a map $T: HP_{even} (D) \ra HC_{2n}(D)$ by 
$$
T:[f] \mapsto [(f_0, \dots, f_{2n})], \qquad f_{2i} \in \bar{\Omega}^{2i}\tilde{D},
$$
for any cycle $ f = \{f_{2n}\}_{n\geq 0}$ in $\bar{\Omega}^{even}\tilde{D}$. 
It is clear that $T$ maps even periodic cycles to cycles in $\Tot_{2n} \bar{\Omega}\tilde{D}$. 

$T$ is surjective,  for let us take a cycle $f = (f_0, \dots , f_{2n})$
in $\Tot_{2n}\bar{\Omega}\tilde{D}$. Having embedded $f$ in $\bar{\Omega}^{even}\tilde{D}$, 
we calculate that
$
\tdb\tdB f_{2n} = -\tdB\tdb f_{2n} = \tdB^2 f_{2n-2}= 0
$ 
so that $g_{2n+1} = \tdB f_{2n} \in \bar{\Omega}^{2n+1}\tilde{D}$ is a cycle 
in the Hochschild complex. Since Hochshild homology vanishes for $2n+1 > N$,
there exists 
$f_{2n+2} \in \bar{\Omega}^{2n+2}\tilde{D}$ such that $\tdb f_{2n+2} = -g_{2n+1}$. 
Then $\tdb f_{2n+2} + \tdB f_{2n} = 0$. Proceeding this way we construct a cycle
$F$ in $\bar{\Omega}^{even}\tilde{D}$ such that $T(F) = f$. 

The map $T$ is also injective. Let $F = \{ f_{2n}\}_{n\geq 0}$ be a cycle in 
$\bar{\Omega}^{even}\tilde{D}$. Then $[T(F)]= 0 $ in $HC_{2n}(D)$ 
if and only if there exists
a chain 
$h = (h_1, \dots, h_{2n+1}) \in \Tot_{2n+1} \bar{\Omega}\tilde{D}$ such that
$(\tdb + \tdB) h = T(F)$. Then 
$$
\tdb(f_{2n+2} - \tdB h_{2n+1}) = \tdb f_{2n+2} + \tdB \tdb h_{2n+1}
= - \tdB f_{2n} + \tdB (f_{2n} + \tdB h_{2n-1}) = 0
$$
and so there exists $h_{2n+3} \in \bar{\Omega}^{2n+3} \tilde{D}$ such that 
$\tdb h_{2n+3} + \tdB h_{2n+1} = f_{2n+2}$. This procedure yields
a chain $H\in \bar{\Omega}^{odd}\tilde{D}$ such that $(\tdb + \tdB) H = F$, and so $[F]=0$
in $HP_{even}(D)$, 
proving that $T$ is injective. Note that $HC_{2n}(D) \simeq HC_{2m}(D)$ for all $n,m$ 
such that $2n, 2m >N$. The argument is the same in the odd case.

Returning to the proof of the theorem, we first use Theorem \ref{HHcts} to 
show that $HH_n(A, \bar{\otimes}) = 0$ for all $n > N$. Then using the above
remark and continuity of $HC$ we  write
$$
\dlim HP_{even}(A_\alpha) = \dlim HC_{2n}(A_\alpha)
= HC_{2n}(A, \bar{\otimes}) = HP_{even }(A,
\bar{\otimes}),
$$
provided $2n > N$. The proof of the odd case is the same.

\mbox{}

{\bf 3. Periodic cyclic homology of $S(GL(n))$} 

Let $F$ be a non-archimedean local field and let $G= GL(n) = GL(n, F)$. 
 Let $K$ be a compact open subgroup of $G$.  Define $S(G//K)$ 
to be all functions $f: G \ra \bc$ which are $K$-bi-invariant and 
rapidly decreasing. Then $S(G//K)$, in its standard seminorm topology,
 is a unital nuclear \frechet\ algebra.  The 
Schwartz algebra $S(G)$ is given by
$
S(G) = \bigcup_K S(G//K)
$
in the inductive limit topology \cite{Silberger}. The algebras $S(G)$, $S(G//K)$ satisfy the 
conditions of Theorem \ref{nuke}.
By Mischenko's theorem \cite{Mischenko}, the Fourier transform determines an isomorphism
of unital \frechet\ algebras: 
$$
S(G//K) \simeq \bigoplus_M\left[C^\infty(F(M:K))\right]^{W(M:K)}
$$
where $F(M:K)\ra E_2(M:K)$ is the complex Hermitian vector bundle of $K$-fixed
vectors in the induced Hilbert bundle $F(M)\ra E_2(M)$. The vector bundle 
$F(M:K)$ is trivialized. One Levi subgroup $M$ is chosen in each $G$-conjugacy
class. 

We are led to the following issue. Let $X$ be a compact smooth manifold (in 
fact a compact torus), $W$ a finite group acting on $X$, and $F\ra X$ a
(trivialized) complex Hermitian vector bundle such that $F$ 
is a $W$-bundle. The group $W$ acts via intertwining operators $a(w:x): F_x \ra
F_{wx}$. For the group $GL(n)$, normalized intertwining operators \cite{Shahidi} 
may be chosen
such that each isotropy subgroup $W_x$ acts trivially in the fibre $F_x$, i.e.
$w\in W_x$ implies $a(w:x) = 1$.

We consider $B = C^\infty(X)^W$, $A= [C^\infty(\mbox{\rm End}\; F)]^W$, 
$E= C^\infty(F)^W$. It is elementary to check that $E$ is an $A-B$-bimodule. 
Then we have a map $\Phi: A \ra End_B (E)$. 
\begin{lemma}\label{Lemma4}
Let $v\in F_x$. Then there exists an invariant section $s$ of $F$ such that 
$s(x) = v$. 
\end{lemma}
\proof Choose a smooth section $t$ such that $t(x) = |W_x|^{-1}v$ and $supp\, t$ does not
contain any point in the orbit $Wx$ except $x$. Now average $t$ by defining 
$s = \sum_{w\in W} wt$.  Then $s$ is an invariant smooth section such that 
$$
s(x) = \sum_{w\in W} a(w:x) t(w^{-1}x)= 
\sum_{w\in W_x} a(w:x) t(x) + \sum_{w\not\in W_x} a(w:x) t(w^{-1}x)
= |W_x|t(x)  = v.
$$

\begin{lemma}\label{Lemma5}
The map 
$
\Phi: A \ra End_B(E)
$
is an isomorphism of \frechet\ algebras.
\end{lemma}
\proof Injectivity of the map $\Phi$ follows from  Lemma \ref{Lemma4}. To prove
surjectivity, define, for $H\in End_A(E)$,
$
f_H(x)(v) = (Hs_v)(x)
$
with $v\in F_x$, $s_v$ an invariant section through $v$. Given that $(Hs_v)(x)
= H(x)s_v(x)$, this definition is independent of the choice of the section
$s_v$. It is elementary to check that $f_H\in A$.

By Lemma \ref{Lemma5},  the algebras $A$ and $B$ are 
Morita equivalent \frechet\ algebras, which implies that they have the same 
Hochschild homology \cite[p.~194]{Loday}. We now have that
$$
HH_*(S(G//K)) = \bigoplus_M HH_*(C^\infty(X(M:K))^{W(M:K)})
$$
Implicit in the proof of 
Lemma 45, p.~344 of \cite{Connes} is the identification of Hochschild 
homology of the algebra $C^\infty(V)$  with the differential forms on $V$. Noting the 
perfect duality between the complexes of forms and currents 
\cite[p.~44-45]{deRham} and using the invariance result in 
\cite[p.~240]{Wass}, the Hochschild homology of 
$W(M:K)$-invariant smooth functions on the smooth manifold $X(M:K)$ 
may be identified with the $W(M:K)$-invariant differential forms on $X(M:K)$. 
Given that $\dim X(M:K) \leq n$
for all compact open subgroups $K$ of $GL(n)$, we have that $HH_p(S(G//K)) = 0$ for all $p >
n$ and all such $K$. Using Theorem \ref{HPcts} we have established the following
result. 
\begin{theorem}\label{Main}
$$
HP_*(S(G), \bar{\otimes}) = \dlim HP_*(S(G//K))
$$
\end{theorem}
We remark that the homology theory on the right is the same as the theory 
$h_*(S(G)) $ of \cite{BHP}.

Each quotient space $X(M:K)/W(M:K)$ creates a disjoint union of compact orbifolds, which 
together form the tempered dual of $GL(n)$.
Each orbifold is the quotient of a compact torus $\bt^k$ by a
product of symmetric  groups and we have $k\leq n$ \cite{Plymen}. By the 
de Rham cohomology of an orbifold $X/W$ we shall mean the $W$-invariant part
of the de Rham cohomology of $X$. 
We first apply the invariance result in \cite[p.~240]{Wass} and then apply 
the fundamental
result of Connes \cite[Ch.~II, Theorem 46]{Connes} to obtain the following theorem. 
\begin{theorem}
The periodic cyclic homology $HP_0(-, \bar{\otimes})$ (resp. $HP_1(-, \bar{\otimes})$) of
the  Schwartz algebra $S(GL(n))$ is isomorphic to the compactly supported 
even (resp. odd) de Rham cohomology of the tempered dual of $GL(n)$.
\end{theorem}

{\footnotesize This research was supported by grants from the LMS and Exeter University 
Research Fund. We would like to thank  Peter Schneider and Freydoon Shahidi for valuable
discussions.}

\end{document}